\documentclass[11pt]{amsart}
\usepackage{amsmath,amssymb}
\newtheorem{theorem}{Theorem}

\newtheorem{corollary}[theorem]{Corollary}

\theoremstyle{definition}
\newtheorem{definition}[theorem]{Definition}
\newtheorem{remark}[theorem]{Remark}

\begin{document}

\title{The isoperimetric inequality}
\author{Simon Brendle and Michael Eichmair}
\address{Columbia University, 2990 Broadway, New York NY 10027, USA}
\address{University of Vienna, Faculty of Mathematics, Oskar-Morgenstern-Platz 1, 1090 Vienna, Austria, ORCID: 0000-0001-7993-9536}
\maketitle

\section{The isoperimetric inequality and the Sobolev inequality} 

\label{introduction}

The isoperimetric problem is one of the oldest and most famous problems in geometry. Its origins date back to the legend of Queen Dido founding the City of Carthage, as told in Virgil's Aeneid.

In two dimensions, the isoperimetric inequality asserts that a disk has the smallest boundary length among all domains in the plane with a given area. 

\begin{theorem}[Isoperimetric inequality in the plane]
\label{isoperimetric.inequality.2D}
Let $E$ be a compact domain in $\mathbb{R}^2$ with smooth boundary. Then 
\[|\partial E| \geq 2 \, \pi^{\frac{1}{2}} \, |E|^{\frac{1}{2}}.\] 
\end{theorem}

Here, $|E|$ denotes the area of $E$ and $|\partial E|$ denotes the length of the boundary $\partial E$. Note that disks achieve equality in the isoperimetric inequality. Indeed, if $E$ is a closed disk of radius $r$ in the plane, then $|E| = \pi r^2$ and $|\partial E| = 2\pi r$. 

Theorem \ref{isoperimetric.inequality.2D} is a special case of a more general inequality which holds in arbitrary dimension.

\begin{theorem}[Isoperimetric inequality in $\mathbb{R}^n$]
\label{isoperimetric.inequality}
Let $E$ be a compact domain in $\mathbb{R}^n$ with smooth boundary. Then 
\[|\partial E| \geq n \, |B_1^n|^{\frac{1}{n}} \, |E|^{\frac{n-1}{n}}.\] 
\end{theorem}

Here, $|E|$ denotes the volume of $E$ and $|\partial E|$ denotes the $(n-1)$-dimensional measure of the boundary $\partial E$. Moreover, $B_1^n = \{x \in \mathbb{R}^n: |x| < 1\}$ denotes the open unit ball in $\mathbb{R}^n$ and $|B_1^n|$ denotes its volume. 

The isoperimetric inequality is sharp on balls. To see this, recall that the volume and boundary area of the unit ball in $\mathbb{R}^n$ are related by $|\partial B_1^n| = n \, |B_1^n|$. Hence, if $E$ is a closed ball of radius $r$, then $|E| = |B_1^n| \, r^n$ and $|\partial E| = |\partial B_1^n| \, r^{n-1} = n \, |B_1^n| \, r^{n-1}$.

Another important inequality related to the isoperimetric inequality is the sharp version of the Sobolev inequality. 

\begin{theorem}[Sobolev inequality on $\mathbb{R}^n$]
\label{sobolev.inequality}
Let $f$ be a smooth function on $\mathbb{R}^n$ with compact support. Then 
\[\int_{\mathbb{R}^n} |\nabla f| \geq n \, |B_1^n|^{\frac{1}{n}} \, \bigg ( \int_{\mathbb{R}^n} |f|^{\frac{n}{n-1}} \bigg )^{\frac{n-1}{n}}.\] 
\end{theorem}

The Sobolev inequality plays a fundamental role in the modern theory of partial differential equations. For a function defined on a ball, we have the following variant of the Sobolev inequality.

\begin{theorem}[Sobolev inequality on a ball]
\label{sobolev.inequality.on.a.ball} 
Let $f$ be a positive smooth function on the closed unit ball $\bar{B}_1^n$. Then 
\[\int_{B_1^n} |\nabla f| + \int_{\partial B_1^n} f \geq n \, |B_1^n|^{\frac{1}{n}} \, \bigg ( \int_{B_1^n} f^{\frac{n}{n-1}} \bigg )^{\frac{n-1}{n}}.\] 
\end{theorem} 

Note that Theorem \ref{sobolev.inequality.on.a.ball} implies the Sobolev inequality on $\mathbb{R}^n$ (Theorem \ref{sobolev.inequality}). To see this, we assume that $f$ is a smooth function on $\mathbb{R}^n$ with compact support. After a suitable rescaling, we may assume that the support of $f$ is contained in the open unit ball $B_1^n$. We then apply Theorem \ref{sobolev.inequality.on.a.ball} to the function $\sqrt{j^{-2}+f^2}$ and send $j \to \infty$. 

Moreover, Theorem \ref{sobolev.inequality} implies the isoperimetric inequality (Theorem \ref{isoperimetric.inequality}). To see this, we assume that $E$ is a compact domain in $\mathbb{R}^n$ with smooth boundary. We then approximate the indicator function of $E$ by a sequence of nonnegative smooth functions with compact support. To explain this, we fix a smooth cutoff function $\eta: [0,\infty) \to [0,\infty)$ such that $\eta(s)=1$ for $s \in [0,1]$, $\eta'(s) \leq 0$ for $s \in [1,2]$, and $\eta(s)=0$ for $s \in [2,\infty)$. For each positive integer $j$, we define $f_j(x) = \eta(j \, \text{\rm dist}(x,E))$. If $j$ is sufficiently large, then $f_j$ is a nonnegative smooth function on $\mathbb{R}^n$. Moreover, 
\[\int_{\mathbb{R}^n} f_j^{\frac{n}{n-1}} \to |E|,\] 
while 
\[\int_{\mathbb{R}^n} |\nabla f_j| \to |\partial E|\] 
as $j \to \infty$. Theorem \ref{sobolev.inequality} then implies $|\partial E| \geq n \, |B_1^n| \, |E|^{\frac{n-1}{n}}$. 

In Sections \ref{transport.proofs} and \ref{abp.proof} we present several different proofs of Theorem \ref{sobolev.inequality.on.a.ball}. In Section \ref{transport.proofs}, we sketch how Theorem \ref{sobolev.inequality.on.a.ball} can be proven using measure transportation. This strategy is due to Gromov and can be implemented in two ways. Gromov's original approach uses the Knothe rearrangement. An alternative approach, due to McCann and Trudinger, is based on the Monge-Amp\`ere equation. In Section \ref{abp.proof}, we discuss a proof of Theorem \ref{sobolev.inequality.on.a.ball} due to Cabr\'e that uses linear partial differential equations and the Alexandrov-Bakelman-Pucci method. 

\section{Proof of Theorem \ref{sobolev.inequality.on.a.ball} using measure transportation} 

\label{transport.proofs}

In this section, we present the measure transportation approach to Theorem \ref{sobolev.inequality.on.a.ball}. By scaling, one can reduce to the special case where $\int_{B_1^n} f^{\frac{n}{n-1}} = |B_1^n|$. The first step of the proof involves constructing a smooth map $\Phi$ from the open unit ball $B_1^n$ into itself with the following properties: 
\begin{itemize}
\item[(i)] For each point $x \in B_1^n$, the eigenvalues of the differential $D\Phi(x)$ are nonnegative real numbers. 
\item[(ii)] For each point $x \in B_1^n$, the determinant $\det D\Phi(x)$ equals $f(x)^{\frac{n}{n-1}}$. 
\end{itemize}
Suppose that $\Phi$ is a map with these properties. We may view $\Phi$ as a vector field defined on $B_1^n$. Since the eigenvalues of the differential $D\Phi$ are nonnegative real numbers, their geometric mean can be estimated from above by their arithmetic mean. This gives 
\[n \, f^{\frac{1}{n-1}} = n \, (\det D\Phi)^{\frac{1}{n}} \leq \text{\rm tr}(D\Phi) = \text{\rm div} \, \Phi\] 
at each point in $B_1^n$. Since $\Phi$ takes values in the unit ball, we know that $-\langle \nabla f,\Phi \rangle \leq |\nabla f|$ at each point in $B_1^n$. Consequently, 
\[n \, f^{\frac{n}{n-1}} \leq f \, \text{\rm div} \, \Phi = \text{\rm div}(f \Phi) - \langle \nabla f,\Phi \rangle \leq \text{\rm div}(f \Phi) + |\nabla f|\] 
at each point in $B_1^n$. In the next step, we integrate over the ball $B_r^n = \{x \in \mathbb{R}^n: |x| < r\}$, where $0 < r < 1$. Using the divergence theorem, we conclude that 
\begin{align*} 
n \int_{B_r^n} f^{\frac{n}{n-1}} 
&\leq \int_{B_r^n} \text{\rm div}(f \Phi) + \int_{B_r^n} |\nabla f| \\ 
&= \int_{\partial B_r^n} f \, \Big \langle \Phi,\frac{x}{r} \Big \rangle + \int_{B_r^n} |\nabla f| 
\end{align*} 
for each $0 < r < 1$. On the other hand, using again the fact that $\Phi$ maps into the unit ball, we obtain $\langle \Phi(x),\frac{x}{r} \rangle \leq 1$ for each $0 < r < 1$ and each point $x \in \partial B_r^n$. This implies 
\[n \int_{B_r^n} f^{\frac{n}{n-1}} \leq \int_{\partial B_r^n} f + \int_{B_r^n} |\nabla f|\] 
for each $0 < r < 1$. Sending $r \to 1$, one obtains 
\[n \int_{B_1^n} f^{\frac{n}{n-1}} \leq \int_{\partial B_1^n} f + \int_{B_1^n} |\nabla f|.\] 
Using the normalization $\int_{B_1^n} f^{\frac{n}{n-1}} = |B_1^n|$, it follows that 
\[n \, |B_1^n|^{\frac{1}{n}} \, \bigg ( \int_{B_1^n} f^{\frac{n}{n-1}} \bigg )^{\frac{n-1}{n}} \leq \int_{\partial B_1^n} f + \int_{B_1^n} |\nabla f|,\] 
as desired.

It remains to construct a map $\Phi$ that satisfies the conditions (i) and (ii) above. Gromov's proof in \cite{Milman-Schechtman} is based on the Knothe rearrangement \cite{Knothe}. This construction gives a smooth map $\Phi$ from the open unit ball $B_1^n$ to itself with the following properties: 
\begin{itemize}
\item For each point $x \in B_1^n$, the differential $D\Phi(x)$ is a triangular matrix and the diagonal entries of $D\Phi(x)$ are nonnegative. 
\item For each point $x \in B_1^n$, the determinant $\det D\Phi(x)$ equals $f(x)^{\frac{n}{n-1}}$. 
\end{itemize} 
Clearly, the Knothe map $\Phi$ satisfies the conditions (i) and (ii) above.

Let us sketch the construction of the Knothe map. For simplicity, we consider the special case $n=2$. The Knothe map $\Phi: B_1^2 \to B_1^2$ has the form $\Phi(x_1,x_2) = (\varphi_1(x_1),\varphi_2(x_1,x_2))$ for $(x_1,x_2) \in B_1^2$. The function $\varphi_1$ maps the interval $(-1,1)$ to itself and satisfies 
\[\frac{\int_{\{(x_1,x_2) \in B_1^2: x_1 \leq s_1\}} f^2}{\int_{B_1^2} f^2} = \frac{|\{(\xi_1,\xi_2) \in B_1^2: \xi_1 \leq \varphi_1(s_1)\}|}{|B_1^2|}\] 
for each $s_1 \in (-1,1)$. For each $s_1 \in (-1,1)$, the function $x_2 \mapsto \varphi_2(s_1,x_2)$ maps the interval $(-\sqrt{1-s_1^2},\sqrt{1-s_1^2})$ to the interval $(-\sqrt{1-\varphi_1(s_1)^2},\sqrt{1-\varphi_1(s_1)^2})$ and satisfies 
\[\frac{\int_{\{(x_1,x_2) \in B_1^2: x_1 = s_1, x_2 \leq s_2\}} f^2}{\int_{\{(x_1,x_2) \in B_1^2: x_1 = s_1\}} f^2} = \frac{|\{(\xi_1,\xi_2) \in B_1^2: \xi_1 = \varphi_1(s_1), \xi_2 \leq \varphi_2(s_1,s_2)\}|}{|\{(\xi_1,\xi_2) \in B_1^2: \xi_1 = \varphi_1(s_1)\}|}\] 
for each $s_2 \in (-\sqrt{1-s_1^2},\sqrt{1-s_1^2})$. 

We next describe an alternative approach, due to McCann and Trudinger, which is based on a different choice of the map $\Phi$. The key step in this approach is to solve a suitable boundary value problem for the Monge-Amp\`ere equation. It was shown by Caffarelli \cite{Caffarelli} and Urbas \cite{Urbas} that there exists a convex function $u: \bar{B}_1^n \to \mathbb{R}$ with the following properties: 
\begin{itemize}
\item The function $u$ is smooth and solves the Monge-Amp\`ere equation 
\[\det D^2 u = f^{\frac{n}{n-1}}\] 
at each point in $\bar{B}_1^n$. 
\item The gradient map 
\[x \mapsto \nabla u(x)\] 
maps $\bar{B}_1^n$ to itself. 
\end{itemize}
We now define $\Phi$ to be the gradient map of $u$, so that $\Phi(x) = \nabla u(x)$ for each $x \in B_1^n$. At each point $x \in B_1^n$, the differential $D\Phi(x)$ is a symmetric matrix with nonnegative eigenvalues, and the determinant $\det D\Phi(x)$ equals $f(x)^{\frac{n}{n-1}}$. Therefore, the gradient map $\Phi$ satisfies the conditions (i) and (ii) above.

\begin{remark} 
The solution of the Monge-Amp\`ere equation has a natural interpretation in terms of optimal mass transport (see \cite{Brenier}, \cite{McCann-Guillen}). To explain this, let $u$ denote the solution of the Monge-Amp\`ere equation described above. Let $\mu$ denote the measure on $\bar{B}_1^n$ which has density $f^{\frac{n}{n-1}}$ with respect to the Lebesgue measure. Let $\nu$ denote the Lebesgue measure on $\bar{B}_1^n$. Note that $\mu(\bar{B}_1^n) = \nu(\bar{B}_1^n)$ in view of our normalization. We then consider the problem of minimizing the transport cost 
\[\frac{1}{2} \int_{\bar{B}_1^n \times \bar{B}_1^n} |x-\xi|^2 \, d\pi(x,\xi)\] 
over all measures $\pi$ on $\bar{B}_1^n \times \bar{B}_1^n$ with the property that the marginal distributions of $\pi$ are given by $\mu$ and $\nu$. It is known that there exists a measure $\pi$ which minimizes the transport cost. Moreover, the optimal measure $\pi$ is supported on the graph $\{(x,\nabla u(x)) \in \bar{B}_1^n \times \bar{B}_1^n: x \in \bar{B}_1^n\}$. 
\end{remark}

\section{Proof of Theorem \ref{sobolev.inequality.on.a.ball} using the Alexandrov-Bakelman-Pucci method} 

\label{abp.proof}

In this section, we describe a proof of Theorem \ref{sobolev.inequality.on.a.ball} using the Alexandrov-Bakelman-Pucci technique. This technique plays a central role in the theory of partial differential equations, where it is used to prove a-priori estimates for elliptic partial differential equations in non-divergence form. Cabr\'e \cite{Cabre} showed that the Alexandrov-Bakelman-Pucci technique can be used to give an alternative proof of the isoperimetric inequality. His argument can be adapted to give a proof of the Sobolev inequality. 

By scaling, one can reduce to the special case where $\int_{B_1^n} |\nabla f| + \int_{\partial B_1^n} f = n \int_{B_1^n} f^{\frac{n}{n-1}}$. This normalization ensures that one can find a function $u: \bar{B}_1^n \to \mathbb{R}$ with the following properties: 
\begin{itemize} 
\item The function $u$ is twice continuously differentiable and solves the linear partial differential equation 
\[\text{\rm div}(f \, \nabla u) = n \, f^{\frac{n}{n-1}} - |\nabla f|\] 
at each point in $\bar{B}_1^n$.
\item The function $u$ satisfies the Neumann boundary condition 
\[\langle \nabla u(x),x \rangle = 1\] 
at each point $x \in \partial B_1^n$. 
\end{itemize}
The existence and regularity of $u$ follow from the standard theory of linear elliptic partial differential equations of second order.

Let $\Phi: B_1^n \to \mathbb{R}^n$ denote the gradient map of $u$, so that $\Phi(x) = \nabla u(x)$ for each $x \in B_1^n$. Let $A$ denote the set of all points $x \in B_1^n$ with the property that $|\nabla u(x)| < 1$ and the Hessian $D^2 u(x)$ is weakly positive definite. 

Clearly, $-\langle \nabla f,\nabla u \rangle \leq |\nabla f|$ at each point in $A$. The partial differential equation for $u$ implies that  
\[f \, \Delta u = \text{\rm div}(f \, \nabla u) - \langle \nabla f,\nabla u \rangle \leq \text{\rm div}(f \, \nabla u) + |\nabla f| = n \, f^{\frac{n}{n-1}}\] 
at each point in $A$. Applying the arithmetic-geometric mean inequality to the eigenvalues of the Hessian of $u$, one obtains 
\[0 \leq \det D^2 u \leq \Big ( \frac{\Delta u}{n} \Big )^n \leq f^{\frac{n}{n-1}}\] 
at each point in $A$. Using the change-of-variables formula, one can estimate the measure of the image $\Phi(A)$. This gives 
\begin{equation} 
\label{upper.bound.for.measure.of.Phi(A)}
|\Phi(A)| \leq \int_A |\det D\Phi| = \int_A |\det D^2 u| \leq \int_A f^{\frac{n}{n-1}} \leq \int_{B_1^n} f^{\frac{n}{n-1}}. 
\end{equation} 
On the other hand, it can be shown that the set $\Phi(A)$ contains the open unit ball $B_1^n$. To see this, suppose that a point $\xi \in B_1^n$ is given. It follows from the Neumann boundary condition for $u$ that the function $x \mapsto u(x) - \langle x,\xi \rangle$ attains its minimum at an interior point $x_0 \in B_1^n$. The first and second order conditions at the minimum point imply that $\nabla u(x_0) = \xi$ and the Hessian $D^2 u(x_0)$ is weakly positive definite. Thus, $x_0 \in A$ and $\Phi(x_0) = \xi$.  

Since $\Phi(A)$ contains the open unit ball $B_1^n$, we obtain 
\begin{equation} 
\label{lower.bound.for.measure.of.Phi(A)}
|\Phi(A)| \geq |B_1^n|. 
\end{equation}
Combining (\ref{upper.bound.for.measure.of.Phi(A)}) and (\ref{lower.bound.for.measure.of.Phi(A)}) gives 
\[\int_{B_1^n} f^{\frac{n}{n-1}} \geq |B_1^n|.\] 
In view of the normalization, it follows that 
\[\int_{B_1^n} |\nabla f| + \int_{\partial B_1^n} f = n \int_{B_1^n} f^{\frac{n}{n-1}} \geq n \, |B_1^n|^{\frac{1}{n}} \, \bigg ( \int_{B_1^n} f^{\frac{n}{n-1}} \bigg )^{\frac{n-1}{n}}.\] 
This completes the proof of Theorem \ref{sobolev.inequality.on.a.ball}. 

\section{The Sobolev inequality and the isoperimetric inequality on a hypersurface in $\mathbb{R}^{n+1}$} 

\label{hypersurfaces}

We next discuss how the Sobolev inequality and the isoperimetric inequality can be generalized to hypersurfaces in $\mathbb{R}^{n+1}$. It is particularly natural to study this question for minimal hypersurfaces. 

To explain the notion of a minimal hypersurface, we first recall the definition of the mean curvature. Suppose that $\Sigma$ is a compact smooth hypersurface in $\mathbb{R}^{n+1}$ (possibly with boundary), and let $p$ be a point on $\Sigma$. We may locally write $\Sigma$ as a level set $w(x_1,\hdots,x_{n+1})=0$, where $w$ is a smooth function which is defined on an open neighborhood of $p$ and satisfies $\nabla w \neq 0$. The unit normal vector field to $\Sigma$ is given by $\nu = \frac{\nabla w}{|\nabla w|}$. Moreover, the mean curvature of $\Sigma$ is given by 
\[H = \frac{\Delta w - (D^2 w)(\nu,\nu)}{|\nabla w|} = \frac{\Delta w}{|\nabla w|} - \frac{(D^2 w)(\nabla w,\nabla w)}{|\nabla w|^3}.\] 
It turns out that this definition depends only on the hypersurface $\Sigma$ and the choice of orientation. It does not, however, depend on the choice of the defining function $w$. 

The notion of mean curvature is closely related to the formula for the first variation of area. To explain this, suppose that $V$ is a smooth vector field on $\mathbb{R}^{n+1}$. If $\Sigma$ has non-empty boundary, we assume that the vector field $V$ vanishes along the boundary of $\Sigma$. We consider the deformed hypersurfaces $\Sigma_s = \varphi_s(\Sigma)$, where $s$ is a small real number and the maps $\varphi_s: \mathbb{R}^{n+1} \to \mathbb{R}^{n+1}$ are defined by $\varphi_s(x) = x+s \, V(x)$ for $x \in \mathbb{R}^{n+1}$. In other words, we deform the hypersurface $\Sigma$ with a velocity given by the vector field $V$. Since $V$ vanishes along the boundary of $\Sigma$, this deformation leaves the boundary of $\Sigma$ unchanged. With this understood, the first order change in the area is given by 
\[\frac{d}{ds} |\Sigma_s| \Big |_{s=0} = \int_\Sigma H \, \langle V,\nu \rangle,\] 
where $H$ denotes the mean curvature of $\Sigma$ (see, e.g., \cite{Colding-Minicozzi}, Chapter 1, Section 1).

\begin{definition} 
We say that $\Sigma$ is a minimal hypersurface if the mean curvature of $\Sigma$ vanishes identically.
\end{definition}

In particular, if $\Sigma$ is a minimal hypersurface, then $\Sigma$ is a critical point of the area functional. 

There are many examples of minimal surfaces in $\mathbb{R}^3$. The most basic ones are the plane, the catenoid 
\[\{(x_1,x_2,x_3) \in \mathbb{R}^3: x_1^2+x_2^2 - \cosh^2(x_3) = 0\},\] 
and the helicoid 
\[\{(x_1,x_2,x_3) \in \mathbb{R}^3: x_1 \sin(x_3) - x_2 \cos(x_3) = 0\}.\] 

In 1921, Carleman \cite{Carleman} showed that the isoperimetric inequality holds for two-dimensional minimal surfaces in $\mathbb{R}^3$ that are diffeomorphic to a disk.

\begin{theorem}[Isoperimetric inequality for minimal disks]
\label{isoperimetric.inequality.for.minimal.disks}
Let $\Sigma$ be a compact two-dimensional minimal surface in $\mathbb{R}^3$ with boundary $\partial \Sigma$. If $\Sigma$ is diffeomorphic to a disk, then 
\[|\partial \Sigma| \geq 2 \, \pi^{\frac{1}{2}} \, |\Sigma|^{\frac{1}{2}}.\] 
\end{theorem}

Note that this inequality is sharp. Carleman's proof of Theorem \ref{isoperimetric.inequality.for.minimal.disks} uses techniques from complex analysis.

Theorem \ref{isoperimetric.inequality.for.minimal.disks} raises the question whether the isoperimetric inequality holds for minimal surfaces of arbitrary dimension and topology. In the 1970s, Allard \cite{Allard} and Michael and Simon \cite{Michael-Simon} proved a general Sobolev inequality which holds for arbitrary hypersurfaces in Euclidean space (and, more generally, for submanifolds of arbitrary codimension). Their arguments are based on the monotonicity formula in minimal surface theory together with covering arguments. More recently, Castillon \cite{Castillon} gave an alternative proof based on techniques from optimal transport. However, these works do not give a sharp constant. In 2019, the first-named author proved a sharp version of the Michael-Simon-Sobolev inequality. 

\begin{theorem}[Sobolev inequality on a hypersurface] 
\label{sharp.michael.simon.inequality}
Let $\Sigma$ be a compact hypersurface in $\mathbb{R}^{n+1}$ with boundary $\partial \Sigma$. Let $f$ be a positive smooth function on $\Sigma$. Then 
\[\int_\Sigma \sqrt{|\nabla^\Sigma f|^2 + f^2 \, H^2} + \int_{\partial \Sigma} f \geq n \, |B_1^n|^{\frac{1}{n}} \, \bigg ( \int_\Sigma f^{\frac{n}{n-1}} \bigg )^{\frac{n-1}{n}}.\] 
Here and below, $\nabla^\Sigma f$ denotes the gradient of $f$ along $\Sigma$. 
\end{theorem}

In particular, if $\Sigma$ is a minimal hypersurface, then the mean curvature term vanishes and we can draw the following conclusion. 

\begin{corollary}[Sobolev inequality on a minimal hypersurface] 
\label{sobolev.inequality.on.a.minimal.hypersurface}
Let $\Sigma$ be a compact minimal hypersurface in $\mathbb{R}^{n+1}$ with boundary $\partial \Sigma$. Let $f$ be a positive smooth function on $\Sigma$. Then 
\[\int_\Sigma |\nabla^\Sigma f| + \int_{\partial \Sigma} f \geq n \, |B_1^n|^{\frac{1}{n}} \, \bigg ( \int_\Sigma f^{\frac{n}{n-1}} \bigg )^{\frac{n-1}{n}}.\] 
\end{corollary}

Putting $f=1$, we obtain the following result.

\begin{corollary}[Isoperimetric inequality on a minimal hypersurface] 
\label{isoperimetric.inequality.on.a.minimal.hypersurface}
Let $\Sigma$ be a compact minimal hypersurface in $\mathbb{R}^{n+1}$ with boundary $\partial \Sigma$. Then 
\[|\partial \Sigma| \geq n \, |B_1^n|^{\frac{1}{n}} \, |\Sigma|^{\frac{n-1}{n}}.\] 
\end{corollary}

Theorem \ref{sharp.michael.simon.inequality} can be proven in two ways. The original proof by the first-named author \cite{Brendle} uses the Alexandrov-Bakelman-Pucci technique. This involves studying a Neumann boundary value problem for a linear partial differential equation on $\Sigma$. The authors \cite{Brendle-Eichmair} recently gave an alternative proof which uses optimal transport. In \cite{Brendle-Eichmair}, only the special case $f=1$ is considered, but the proof can be adapted so that it works for an arbitrary positive smooth function $f$.

In the optimal transport approach, it is convenient to normalize $f$ such that $\int_\Sigma f^{\frac{n}{n-1}} = 1$. Let $\mu$ denote the measure on $\Sigma$ which has density $f^{\frac{n}{n-1}}$ with respect to the volume measure on $\Sigma$. Let $\rho: [0,\infty) \to (0,\infty)$ be a continuous function with the property that $\int_{\bar{B}_1^{n+1}} \rho(|\xi|^2) \, d\xi = 1$. Let $\nu$ denote the measure on the $(n+1)$-dimensional unit ball $\bar{B}_1^{n+1}$ which has density $\rho(|\xi|^2)$ with respect to the $(n+1)$-dimensional Lebesgue measure. By definition, $\mu$ is a probability measure on $\Sigma$ and $\nu$ is a probability measure on $\bar{B}_1^{n+1}$.  

The key idea is to consider the optimal transport problem between $(\Sigma,\mu)$ and $(\bar{B}_1^{n+1},\nu)$, with a quadratic cost function. In other words, we minimize the transport cost 
\[\frac{1}{2} \int_{\Sigma \times \bar{B}_1^{n+1}} |x-\xi|^2 \, d\pi(x,\xi)\] 
over all measures $\pi$ on $\Sigma \times \bar{B}_1^{n+1}$ with the property that the marginal distributions of $\pi$ are given by $\mu$ and $\nu$. Note that this is a transport problem between spaces of different dimensions. 

The solution of the optimal transport problem can be described in terms of a function $u: \Sigma \to \mathbb{R}$. The function $u$ is Lipschitz continuous with Lipschitz constant $1$. Moreover, $u$ is the restriction to $\Sigma$ of a convex function on $\mathbb{R}^{n+1}$. In particular, $u$ is a semi-convex function on $\Sigma$. A classical theorem of Alexandrov implies that $u$ admits first and second derivatives at almost every point on $\Sigma$. We next establish a pointwise inequality involving the Alexandrov Hessian of $u$ and the second fundamental form of $\Sigma$ (compare \cite{Brendle-Eichmair}, Proposition 9). This inequality can viewed as the analogue of the Monge-Amp\`ere equation in this setting. Using the arithmetic-geometric mean inequality, we obtain a pointwise inequality involving the Laplacian of $u$ and the mean curvature of $\Sigma$ (compare \cite{Brendle-Eichmair}, Corollary 10). More precisely, we can show that the inequality 
\begin{equation} 
\label{Laplacian}
n \, \alpha^{-\frac{1}{n}} \, f^{\frac{n}{n-1}} \leq f \, \Delta_\Sigma u + \langle \nabla^\Sigma f,\nabla^\Sigma u \rangle + \sqrt{|\nabla^\Sigma f|^2 + f^2 \, H^2} 
\end{equation} 
holds almost everywhere on $\Sigma$. Here, $\Delta_\Sigma u$ denotes the trace of the Alexandrov Hessian of $u$. Moreover, $\alpha$ is defined by 
\[\alpha = \sup_{z \in [0,1)} \int_{-\sqrt{1-z^2}}^{\sqrt{1-z^2}} \rho(z^2+y^2) \, dy.\] 
Note that $\alpha$ is a positive real number that depends on our choice of the density $\rho$. Integrating the inequality (\ref{Laplacian}) over $\Sigma$ gives 
\[n \, \alpha^{-\frac{1}{n}} \int_\Sigma f^{\frac{n}{n-1}} \leq \int_{\partial \Sigma} f + \int_\Sigma \sqrt{|\nabla^\Sigma f|^2 + f^2 \, H^2}.\]
Using the normalization $\int_\Sigma f^{\frac{n}{n-1}} = 1$, it follows that 
\begin{equation} 
\label{inequality.involving.alpha}
n \, \alpha^{-\frac{1}{n}} \, \bigg ( \int_\Sigma f^{\frac{n}{n-1}} \bigg )^{\frac{n-1}{n}} \leq \int_{\partial \Sigma} f + \int_\Sigma \sqrt{|\nabla^\Sigma f|^2 + f^2 \, H^2}. 
\end{equation} 
Finally, one needs to make a suitable choice of the density $\rho$. For each positive integer $j$, we define a continuous density $\rho_j: [0,\infty) \to (0,\infty)$ by 
\[\rho_j(s) = \frac{1}{c_j \sqrt{\max \{1-s,j^{-1}\}}}\]
for all $s \geq 0$, where 
\[c_j = \int_{\bar{B}_1^{n+1}} \frac{1}{\sqrt{\max \{1-|\xi|^2,j^{-1}\}}} \, d\xi.\] 
This choice of the constant $c_j$ ensures that $\int_{\bar{B}_1^{n+1}} \rho_j(|\xi|^2) \, d\xi = 1$. Note that  
\[\lim_{j \to \infty} c_j = \int_{B_1^{n+1}} \frac{1}{\sqrt{1-|\xi|^2}} \, d\xi = \pi \, |B_1^n|.\] 
Moreover, if we put 
\[\alpha_j = \sup_{z \in [0,1)} \int_{-\sqrt{1-z^2}}^{\sqrt{1-z^2}} \rho_j(z^2+y^2) \, dy,\] 
then 
\[\alpha_j \leq \sup_{z \in [0,1)} \int_{-\sqrt{1-z^2}}^{\sqrt{1-z^2}} \frac{1}{c_j \sqrt{1-z^2-y^2}} \, dy = \frac{\pi}{c_j}\] 
for each $j$. Consequently, 
\begin{equation} 
\label{limit.of.alpha_j}
\limsup_{j \to \infty} \alpha_j \leq \frac{1}{|B_1^n|}. 
\end{equation} 
Theorem \ref{sharp.michael.simon.inequality} follows by combining (\ref{inequality.involving.alpha}) and (\ref{limit.of.alpha_j}).

\section{Outlook: Isoperimetric problems in Riemannian geometry and mathematical relativity} 

There is an extensive literature concerning isoperimetric problems in Riemannian manifolds; in particular, there is a version of the isoperimetric inequality in hyperbolic space and in the standard sphere (see, e.g., \cite{Burago-Zalgaller}). Gromov proved an isoperimetric inequality which holds for every Riemannian manifold with Ricci curvature bounded below by a positive constant (see \cite{Gromov}, Appendix C). Klartag gave an alternative proof of Gromov's isoperimetric inequality (see \cite{Klartag}, Proposition 5.4). His approach uses needle decompositions; to construct these, one considers a solution of an optimal transport problem, where the cost function is given by the Riemannian distance. Bray \cite{Bray} used isoperimetric surfaces to prove volume comparison theorems for three-dimensional manifolds with lower bounds on the scalar curvature and the Ricci curvature. Isoperimetric surfaces have also found important applications in mathematical general relativity, where they have been shown to mediate between positive energy density on small scales and positive mass at infinity in initial data of the Einstein field equations that model an isolated gravitational system (see \cite {Bray}, \cite{Eichmair-Metzger}, \cite{Chodosh-Eichmair-Shi-Yu}).

\end{document}